\documentclass[preprint,10pt]{article}
\usepackage{cite}
\usepackage{url}
\usepackage{ifthen}
\usepackage{multicol}
\usepackage{amsmath,amssymb,amsthm}
\usepackage{latexsym,amssymb,amsmath}
\usepackage[utf8]{inputenc}
\usepackage{graphicx}
\usepackage{epstopdf}
\usepackage{subfig}
\usepackage{float}
\usepackage{color}
\newtheorem{theorem}{Theorem}
\newtheorem{corollary}{Corollary}
\newtheorem{lemma}{Lemma}
\newtheorem{definition}{Definition}

\newtheorem{remark}{Remark}

\numberwithin{equation}{section}
\newboolean{publ}

\newcommand{\ncm}{\newcommand}
\ncm{\maint}{\rlap{$\,\hspace{1.5790pt}=$}\int}
\ncm{\mcint}{\rlap{$\,\hspace{1.5790pt}-$}\int}

\title{On a fractional class of analytic function defined by using a new operator }
\author{  Zainab E. Abdulnaby$^{1}$,  Rabha W. Ibrahim$^{2}$ and  Adem Kilicman$^{1}$   \\
          \small $^1$ Department of Mathematics, Universiti Putra Malaysia, \\
         \small  43400 UPM Serdang, Selangor, Malaysia\\
          \small $^2$ Faculty of Computer Science and Information Technology\\
          \small University Malaya, 50603, Malaysia
}
\begin{document}
\maketitle
\abstract{\noindent In this article, we impose a new class of
fractional analytic functions in the open unit disk. By considering
this class, we define a fractional operator, which is generalized
Salagean and Ruscheweyh differential operators. Moreover, by means
of this
 operator, we introduce an interesting subclass of functions which are analytic and univalent.
  Furthermore, this effort covers coefficient bounds, distortions theorem, radii of starlikeness,
   convexity, bounded turning, extreme points and integral means inequalities of functions
    belongs to this class. Finally, applications involving certain fractional operators are illustrated.

 \section{Introduction and preliminaries}
 Recently, One of the substantive issues in many applications of geometric function theory is how
  to employ the fractional operators to analytic and univalent functions and what the advantages
  for this utilized. In other hand, make use of functional analytic functions to define fractional
   operators and what of the results of this. So far, many mathematicians in different stages
   considered this issues and gave numerous applications based on certain fractional operators
    of analytic function in physics, engineering and mathematical applications (see
    \cite{z13}).

    \bigskip \noindent In the theory of geometric functions, the Koebe function is formulated by

\[f(z)=\frac{z}{(1 - z)^2}=\sum_{n=1}^\infty n z^n. \]
The rotated Koebe function is

\[f_\alpha(z)=\frac{z}{(1-\alpha z)^2}=\sum_{n=1}^\infty
n\alpha^{n-1} z^n\] with a complex number of absolute value 1. The
Koebe function and its rotations are schlicht: that is, univalent
(analytic and one-to-one) and achieving $f(0) = 0$ and $f'(0) = 1.$
Srivastava et al \cite{z6}, introduced a fractional analytic
function as follows:

\[f(z)= \frac{z^{\alpha+1}}{(1-z)^\alpha}, \quad \alpha \in \mathbb{R}.\]
In this effort, we define a class $ \mathcal{A}_{\mu}$ of functional
fractional analytic functions $ F_\mu(z)$ in unit
  disk $\mathbb{U}:= \left\lbrace z \in \mathbb{C};\, |z| < 1\right\rbrace$  as follows:

 \begin{equation}
 F_\mu(z)= \frac{z^\mu}{1-z^\mu},
 \end{equation} where $\mu:= \frac{n+m-1}{m}, \, n,m \in \mathbb{N}.$
Hence, $\mu=1,$ when $n=1$ and has the power series formal:

 \begin{align}\label{fu2}
  F_\mu(z)= z + \sum_{n=2}^{\infty} a_{n} z^{\mu n }   \end{align}
    $$(\mu \geq 1; n \in N, \, z \in \mathbb{U}),$$ which normalize
  by $ F_\mu(z)|_{z=0}=1$ and $ F_\mu^{\prime}(z)|_{z=0}=1$ for
   all $ z \in \mathbb{U}$.

   \bigskip \noindent Recall that a
    function $ F_\mu  \in \mathcal{A}_{\mu}$ is called bounded turning if it
    satisfies the following inequality:

  \begin{align}\label{S}
  \Re \lbrace F_{\mu}^{\prime}(z) \rbrace > \psi \quad (0 \leq \psi < 1), \end{align} and a function
  $ F_{\mu} \in \mathcal{A}_{\mu}$ is starlike function in $\mathbb{U}$ if satisfies

       \begin{align}\label{S}
  \Re \lbrace \frac{ z F_{\mu}^{\prime}(z)}{F_{\mu}(z)} \rbrace > \psi \quad (0 \leq \psi < 1). \end{align}
  Furthermore, a function  $ F_{\mu} \in \mathcal{A}_{\mu}$ is convex function in $ \mathbb{U}$ if satisfies
     \begin{align}\label{S}
  \Re \lbrace 1+ \frac{ z  F_{\mu}^{\prime \prime}(z) }{F_{\mu}^{\prime}(z)} \rbrace >
  \psi \quad (0 \leq \psi < 1). \end{align} (see, for more details \cite{z12} and \cite{z8}).

   \medskip \noindent  Next, if the function $ F_{\mu}(z) $ of form \eqref{fu2}  and $ G _{\mu}(z)= z+ \sum_{n=2}^{\infty}
    b_{n} z^{\mu n}$ are two functions in class $ \mathcal{A}_{\mu}$, then
the convolution (or Hadamard product) of two analytic functions is denoted by $ F_{\mu} * G_{\mu}$ and is given by
    $$ F_{\mu}(z)*G_{\mu}(z)= z+ \sum_{n=2}^{\infty} a_{n} b_{n} z^{\mu n} $$
    and satisfy
   $$ [F_{\mu}(z)*G_{\mu}(z)]\big|_{z=0}= 0 \quad \text{ and }\quad  [(F_{\mu}(z)*G_{\mu}(z))^{\prime}]\big|_{z=0}= 1.$$

\noindent Now let a functional function $ \Theta_{\mu}(z)$  defined as  follows:
    \begin{align}
   \Theta_{\mu}(z) &= \frac{\mu z^{\mu}}{(1-z^{\mu})^{2}} + \frac{\mu z^{\mu}}{(1-z^{\mu})},\nonumber \\
   &= z+ \sum_{n=2}^{\infty} (\mu n) z^{\mu n }, \quad (z \in \mathbb{U}).
   \end{align}
   \noindent   By employing method of the convolution product of analytic function $ \Theta_{\mu}(z)$, we obtain 
   \begin{align}
  \Theta_{\mu,k}(z)& = \underbrace{F_{\mu}(z)*\cdots * F_{\mu}(z)}_{k-times} , \nonumber \\
  & = z+ \sum_{n=2}^{\infty} (\mu n)^{k}  z^{\mu n }, \quad (z \in \mathbb{U}).
\end{align}   

\bigskip \noindent For $ F_{\mu}(z) \in \mathcal{A}_{\mu}$, we define   differential operator $ D_{\mu}^{k}F_{\mu}(z)$ as the following 
 $$ D_{\mu}^{k}F_{\mu}(z)= \Theta_{\mu,k}(z)*F_{\mu}(z), \quad (|z|<1).$$
 where
 \begin{align*}\label{1}
  D_{\mu}^{0}F_{\mu}(z) & = F_{\mu}(z),
  \\
  D_{\mu}^{1}F_{\mu}(z)& =  z F_{\mu}^{\prime}(z) = z+ \sum_{n=2}^{\infty} (\mu n) a_{n} z^{\mu n }  \\
  D_{\mu}^{2}F_{\mu}(z)& =  D \big( D_{\mu} F_{\mu}(z)\big) = z+ \sum_{n=2}^{\infty} (\mu n)^{2} a_{n} z^{\mu n }  \\
  & \qquad \quad \quad \vdots\\
   D_{\mu}^{k}F_{\mu}(z)&= D\left(  D_{\mu}^{k-1}F_{\mu}(z) \right) =z+ \sum_{n=2}^{\infty}  (\mu n)^{k} a_{n} z^{\mu n }.
  \end{align*}
  In general, we write
    \begin{equation}\label{Dm}
    D^{\beta}_{\mu} F_{\mu}(z)=z+ \sum_{n=2}^{\infty}  (\mu n)^{\beta} a_{n}z^{\mu n },
\end{equation}
      $$(\mu\geq 1; n \in \mathbb{N}\setminus \lbrace 1 \rbrace; \beta \in \mathbb{N}_{0}; z \in \mathbb{U}).$$
     \noindent Based of above functional function, we deal with a new operator  $ D_{\mu}^{\beta}F_{\mu}
     (z)$  which is a generalization of the well known
     operators such as Ruscheweyh differential operators in $ \mathbb{U}$ (see \cite{z1}).

 \bigskip
\noindent With  recalling  the principle of subordination between two analytic functions $ f$ and $ g$ in open unit disk $ \mathbb{U}$
 (see\cite{z3} and \cite{z4}), the function  $ f$  is \textit{subordinate} to $g$  if there exists a
  Schwarz function $ w(z)$, analytic in $ \mathbb{U}$  with $ w(z)|_{z=0}=0$ and $ |w(z)|=1$ for all $ z \in \mathbb{U}$ such that
  $$ f(z)= g(w(z)), \quad (z \in \mathbb{U}).$$
In addition, we write this subordination by

  \begin{align}\label{Sub}
 f(z) \prec g(z) \quad (z \in \mathbb{U}).
\end{align}
  \noindent In special case, if the function $ g$ is univalent in $ \mathbb{U}$ then the above subordination is equivalent to
 $$ f(z)|_{z=0}=g(z)|_{z=0} \quad \text{and} \quad f(\mathbb{U}) \subset g(\mathbb{U}).$$
 \noindent
  Now we define the following a new class of analytic functions and investigate
several interesting results.

\begin{definition}\rm{
 Let the functions
  $$ \Phi_\mu(z)= z+ \sum_{n=2}^{\infty} \vartheta_{n} z^{\mu n} \quad \text{and} \quad \Psi_\mu(z)= z+ \sum_{n=2}^{\infty} \lambda_{n} z^{\mu n}$$
  be analytic in the open unit disk $ \mathbb{U}$ where
 $$ \vartheta_{n} \geq 0, \quad \lambda_{n} \geq 0 \quad \text{ and } \quad \vartheta_{n}\geq \lambda_{n}
 \quad ( n \in \mathbb{N}\setminus \left\lbrace0,1\right\rbrace).$$
\noindent Then  a function $ F_\mu(z) \in \mathcal{A}_{\mu}$ is said
to be in the class $
\mathcal{E}_{k,m}(\Phi_\mu,\Psi_\mu,A,B,\mu,\gamma)$ if and only if
$$ \frac{D_\mu^{k}(F_\mu*\Phi_\mu)(z)}{D_\mu^{m}(F_\mu*\Psi_\mu)(z)} \prec (1-\gamma)\frac{1+Az}{1+Bz} +\gamma \quad ( z \in \mathbb{U}),$$
where $ \prec$ represents the subordination in \eqref{Sub}, $
F_\mu(z)*\Psi_\mu(z)\neq 0$, $ A$ and $ B$ are arbitrarily fixed
numbers such that $$ -1 \leq B < A \leq 1\quad \text{and} \quad -1
\leq B < 0$$ with
$$  0\leq \gamma < 1 \quad \text{and} \quad k \geq m \quad ( k, m \in \mathbb{N}_{0}).$$
In other words, $F_\mu\in
\mathcal{E}_{k,m}(\Phi_\mu,\Psi_\mu,A,B,\mu,\gamma)$ if and only if
there exists an analytic function $w(z)$ satisfying
$$w(z)|_{z=0}=0 \quad \text{and} \quad| w(z)|=1 \quad (z \in \mathbb{U})$$
such that

\begin{align}\label{2}
\frac{D_\mu^{k}(F_\mu*\Phi_\mu)(z)}{D_\mu^{m}(F_\mu*\Psi_\mu)(z)}
=(1-\gamma)\frac{1+A\omega(z)}{1+B\omega(z)} +\gamma \quad ( z \in
\mathbb{U}).
\end{align}
The condition  \eqref{2} can be expressed by the equivalent inequality

\begin{align}
\left| \frac{\frac{D_{\mu}^{k}(F_\mu
*\Phi_\mu)(z)}{D_{\mu}^{m}(F_\mu*\Psi_\mu)(z)}-1}
{(A-B)(1-\gamma)- B \left(
\frac{D_{\mu}^{k}(F_\mu*\Phi_\mu)(z)}{D_{\mu}^{m}(F_{\mu}*\Psi_\mu)(z)}-1\right)}
\right| < 1\quad ( z\in \mathbb{U}).
\end{align}
}\end{definition}

 Let $ \mathcal{X}_{\mu}$ be the class of analytic functions $ F_\mu(z)$  in unit disk $ \mathbb{U}$ of the following form
     \begin{equation}\label{-}
    F_\mu(z)= z - \sum_{n=2}^{\infty} a_{n} z^{\mu n },
    \end{equation}
    $$(a_{n}\geq 0; \mu\geq 1; n \in N\setminus \lbrace 1 \rbrace)$$
    which satisfies $ F_\mu(z)|_{z=0 }=0$ and $
F^{\prime}_\mu(z)|_{z=0}=1 $ for all $ z \in \mathbb{U}$. We denote
by $\widetilde{ \mathcal{E}}_{k,m}(\Phi_\mu,\Psi_\mu,A,B,\gamma,\mu)
$ the subclass of functions in  ${
\mathcal{E}}_{k,m}(\Phi_\mu,\Psi_\mu,A,B,\gamma,\mu) $ that has
their non-zero coefficients, from second onwards, all negative.
Further let
$$ \widetilde{ \mathcal{E}}_{k,m}(\Phi_\mu,\Psi_\mu,A,B,\mu,\gamma)
 = { \mathcal{E}}_{k,m}(\Phi_\mu,\Psi_\mu,A,B,\mu,\gamma)\cap \mathcal{X}_{\mu}. $$
For suitable choices of $ \Phi$ and $\Psi$, we definitely obtain the function subclasses of $\mathcal{A}_{\mu}$. For example, we have the following:
\begin{align}\label{8}
 \widetilde{ \mathcal{E}}_{0,0}(\frac{\mu z^{\mu}}{(1-z^{\mu})^{2}}, \frac{z^{\mu}}{1- z^{\mu}},1,-1, \mu,\gamma) = \mathcal{S}_{\mu}^{*}(\gamma)
\end{align} 
and 
\begin{align}\label{9}
\widetilde{ \mathcal{E}}_{0,0}(\frac{\mu^{2} z^{\mu}+ \mu^{2} z^{2\mu}}{(1-z^{\mu})^{3}},\frac{\mu z^{\mu}}{(1-z^{\mu})^{2}}, 1,-1,\mu,\gamma) = \mathcal{K}_{\mu}(\gamma)
\end{align} 
If $ \mu =1$ in \eqref{8} and \eqref{9}, respectively we obtain the well known subclasses:
 $$ \mathcal{S}^{*}(\gamma) \quad \text{and } \quad \mathcal{K}(\gamma),$$
which were investigated by \cite{z9}.

\section{Characterization properties}
In this section, we consider several properties for the function $
F_\mu(z) \in \widetilde{
\mathcal{E}}_{k,m}(\Phi_\mu,\Psi_\mu,A,B,\mu,\gamma)$. We will
divide this section into five subsections.

    \subsection{Coefficient Bounds}
  \begin{theorem}\label{Th1}\rm{
 If $ F_\mu(z) \in \mathcal{A}_{\mu}$ satisfies the following inequality
 \begin{align}\label{Co1}
  \sum_{n=2}^{\infty} \big[(1-B) \left((\mu n)^{k} \vartheta_{n}-(\mu n)^{m}
  \lambda_{n}\right)+(A-B)(1-\gamma)(\mu n)^{m} \lambda_{n}\big]\,| a_{n}|\leq \,(A-B)(1-\gamma).
 \end{align}
 $$( \vartheta_{n} \geq 0;\,\lambda_{n}\geq 0;\,\vartheta_{n}\geq \lambda_{n}\,
 ( n \in \mathbb{N}\setminus{\lbrace 1\rbrace});\,\mu \geq 1;\, 0 \leq \gamma < 1;\,k \geq m;\,k, m \in \mathbb{N}_{0}).$$
  Then $$ F_\mu(z) \in \mathcal{E}_{k,m}(\Phi_\mu,\Psi_\mu,A,B,\mu,\gamma).$$

 \begin{proof}
 Let the condition \eqref{Co1} holds, then we obtain
 \begin{align}\label{Eq2.2}
 & \left| D_{\mu}^{k}(F_\mu*\Phi_\mu)(z)- D_{\mu}^{m}(F_\mu*\Psi_\mu)(z)\right| \nonumber
  \\  &
 \quad \quad  - \left| (A-B)(1-\gamma)D_{\mu}^{m}(F_\mu*\Psi_\mu)(z)- B \left( D_{\mu}^{k}(F_\mu
 *\Phi_\mu)(z) - D_{\mu}^{m}(F_\mu*\Psi_\mu)(z)\right)\right| \nonumber
  \\
 & = \left|\sum_{n=2}^{\infty} \left((\mu n)^{k} \vartheta_{n} - (\mu n)^{m}
 \lambda_{n}\right) a_{n} z^{\mu n} \right|-\bigg| (A-B)(1-\gamma) z +(A-B)
 \nonumber
 \\  &\qquad
 \times(1-\gamma)\sum_{n=2}^{\infty} (\mu n)^{m}
  \vartheta_{n}\, a_{n} z^{\mu n}- B \sum_{n=2}^
  {\infty}\left((\mu n)^{k} \vartheta_{n} - (\mu n)^{m} \lambda_{n}\right)\, a_{n} z^{\mu n} \bigg| \nonumber
  \\  &\leq \sum_{n=2}^{\infty} \left((\mu n)^{k}
   \vartheta_{n} - (\mu n)^{m} \lambda_{n}\right)\,|a_{n}|  r^{\mu n}  +  (A-B)(1-\gamma)\,r \nonumber \\
 & \qquad +  (A-B)(1-\gamma)\sum_{n=2}^{\infty} (\mu n)^{m} \vartheta_{n} |a_{n}| r^{\mu n} + |B | \sum_{n=2}^{\infty}
 \left((\mu n)^{k} \vartheta_{n} - (\mu n)^{m} \lambda_{n}\right)\,|a_{n}|.r^{\mu n}
  \\  &\leq
  \sum_{n=2}^{\infty} \big[(1-B)\left((\mu n)^{k} \vartheta_{n} - (\mu n)^{m} \lambda_{n}\right) \nonumber\\
 & \qquad  + (A-B)(1-\gamma)(\mu n)^{m} \vartheta_{n}\big]\,|a_{n}|-(A-B)(1-\gamma)\,
 \leq\, 0.\label{E2.3}
 \end{align}
 For all $ r( 0 \leq r < 1)$ the inequality  in \eqref{Eq2.2} holds true.  Thus, letting $ r \rightarrow 1-$ in \eqref{Eq2.2}, we obtain \eqref{E2.3}
 Hence,  this completes proof of Theorem  \ref{Th1}.
 \end{proof}
 }\end{theorem}
   \begin{theorem}\label{Th2}\rm{
 If $ F_\mu(z) \in \mathcal{X}_{\mu}$ satisfies the following inequality
 \begin{flalign}\label{Co2}
  \sum_{n=2}^{\infty}\big[(1-B)\left((\mu n)^{k} \vartheta_{n}-(\mu n)^{m}\lambda_{n}\right)
+(A-B)(1-\gamma)(\mu n)^{m}\vartheta_{n}\big]\, a_{n}\leq (A-B)(1-\gamma).
&& \end{flalign}
$$( \vartheta_{n} \geq 0;\,\lambda_{n}\geq 0;\,\vartheta_{n}\geq \lambda_{n}\,
( n \in \mathbb{N}\setminus{\lbrace 1\rbrace});\,\mu \geq 1;\, 0 \leq \gamma < 1;\,k\geq m;\,k, m \in \mathbb{N}_{0}).$$
  Then $$ F_\mu(z) \in \widetilde{ \mathcal{E}}_{k,m}(\Phi_\mu,\Psi_\mu,A,B,\mu,\gamma).$$
 \begin{proof}
 Since
 $$ \widetilde{ \mathcal{E}}_{k,m}(\Phi_\mu,\Psi_\mu,A,B,\mu,\gamma)\subset { \mathcal{E}}_{k,m}(\Phi_\mu,\Psi_\mu,A,B,\mu,\gamma), $$
 we only need to prove the \textit{only if } part of Theorem \ref{Th2} for function $ F_\mu(z) \in \mathcal{X}_{\mu}$ we can write

  \begin{align}
& \left| \frac{\left({D_{\mu}^{k}(F_\mu
*\Phi_\mu)(z)})/({D_{\mu}^{m}(F_\mu*\Psi_\mu)(z)}\right)-1}
{(A-B)(1-\gamma)- B
 \left( ({D_{\mu}^{k}(F_\mu*\Phi_\mu)(z)})/({D_{\mu}^{m}(F_\mu*\Psi_\mu)(z)})-1\right)}  \right| \nonumber
 \\
& = \left| \frac{{D_{\mu}^{k}(F_\mu
*\Phi_\mu)(z)}-{D_{\mu}^{m}(F_\mu*\Psi_\mu)(z)}}
{(A-B)(1-\gamma)D_{\mu}^{m}(F_\mu*\Psi_\mu)(z)- B \left[
{D_{\mu}^{k}(F_\mu*\Phi_\mu)(z)}-{D_{\mu}^{m}(F_\mu*\Psi_\mu)(z)}\right]}
\right| \nonumber
\\
&= \left|\frac{\sum_{n=2}^{\infty}((\mu n)^{k} \vartheta_{n} -
( \mu n)^{m} \lambda_{n} )\, a_{n} z^{\mu n-1}}{  \begin{array}{cc}(A-B)
(1-\gamma)-(A-B)(1-\gamma)\sum_{n=2}^{\infty}(\mu n)^{m} \lambda_{n} a_{n}z^{\mu n -1}
\\
+ B \sum_{n=2}^{\infty}\left((\mu n)^{k} \vartheta_{n} - ( \mu n)^{m} \lambda_{n} \right)
\, a_{n} z^{\mu n-1} \end{array}}   \right| < 1 \nonumber .
\end{align}
Since $ \Re(z) \leq |z|$ for all $ z \in \mathbb{U}$, we have
$$\Re \left\lbrace  \frac{\sum_{n=2}^{\infty}\left((\mu n)^{k} \vartheta_{n} -
( \mu n)^{m} \lambda_{n}\right)\, a_{n} z^{\mu n-1}}{  \begin{array}{cc}(A-B)
(1-\gamma)-(A-B)(1-\gamma)\sum_{n=2}^{\infty}(\mu n)^{m} \lambda_{n} a_{n}z^{\mu n-1}
\\
+ B \sum_{n=2}^{\infty}\left((\mu n)^{k} \vartheta_{n} - ( \mu n)^{m} \lambda_{n}\right)
\, a_{n} z^{\mu n-1} \end{array}}   \right\rbrace   < 1. $$
If we choose $ z $ to be real and let $ z \rightarrow 1-$, we obtain
\begin{align*}
  \sum_{n=2}^{\infty} \big[(1-B)\left((\mu n)^{k} \vartheta_{n}-(\mu n)^{m} \lambda_{n}\right)+(A-B)(1-\gamma)&
  (\mu n)^{m}\vartheta_{n}\big]\, a_{n}\leq (A-B)(1-\gamma).
 \end{align*}
  which is equivalent to \eqref{Co2}.
  The result is sharp for functions $F$ given by \begin{align}\label{fun}
 F_\mu(z)= z -\frac{(A-B)(1-\gamma) }{\big[(1-B) \left((\mu n)^{k} \vartheta_{n}-(\mu n)^{m} \lambda_{n}\right)
 + (A-B)(1-\gamma)(\mu n)^{m} \vartheta_{n}\big]}z^{\mu n}\, ( n \geq 2).
\end{align}
This completes the proof of Theorem \ref{Th2}.
\end{proof}
  } \end{theorem}

\begin{corollary}\rm{
Let a function $F_\mu(z)$ defined by \eqref{-} belongs to
$\widetilde{\mathcal{E}}_{k,m}(\Phi_\mu,\Psi_\mu,A,B,\mu,\gamma) $.
Then
\begin{align}
a_{n} \leq \frac{(A-B)(1-\gamma)}{\big[(1-B) \left((\mu n)^{k} \vartheta_{n} - (\mu n)^{m} \lambda_{n}\right)
 + (A-B)(1-\gamma)(\mu n)^{m} \vartheta_{n}\big]}\quad  (n\geq 2).
\end{align}
}\end{corollary}
 \begin{remark} \rm{By taking different choices for the functions $\Phi_\mu(z)$ and $\Psi_\mu(z)$
 same as stated in \eqref{8} and \eqref{9},
 Theorem \ref{Th2} leads us to the necessary and sufficient conditions for a function $F_\mu(z)$
  to be in the following classes:

$$\mathcal{S}^{*}_{\mu}(\gamma)\quad \text{and} \quad \mathcal{K}_{\mu}(\gamma).$$
}\end{remark}


 \subsection{Distortion theorems}
\begin{theorem}\label{Th:7}  \rm{
Let the function $F_\mu(z)$ defined by \eqref{-} be in
$\widetilde{\mathcal{E}}_{k,m}(\Phi_\mu,\Psi_\mu,A,B,\mu,\gamma)$.
Then

\begin{align}
|F_\mu(z)| & \geq |z|- \frac{(A-B)(1-\gamma)}{\big[(1-B) \left((
2\mu )^{k} \vartheta_{2} - (2\mu )^{m} \lambda_{2}\right)
 + (A-B)(1-\gamma)(2\mu)^{m} \vartheta_{2}\big]}   |z|^{ 2 \mu } \nonumber
 \\
 &  \leq |z|+ \frac{(A-B)(1-\gamma)}{\big[(1-B) \left(( 2\mu )^{k} \vartheta_{2} - (2\mu )^{m} \lambda_{2}\right)
 + (A-B)(1-\gamma)(2\mu)^{m} \vartheta_{2}\big]}   |z|^{ 2 \mu }.
\end{align}
The result is sharp.
\begin{proof}
By considering Theorem \ref{Th1}, since
\begin{align*}
\Xi(n)= \big[(1-B)\left((\mu n)^{k} \vartheta_{n} - (\mu n)^{m} \lambda_{n}\right) + (A-B)(1-\gamma)(\mu n)^{m} \vartheta_{n}\big]
\end{align*}
is an increasing function of $ n \,(n \geq 2)$, we get
\begin{align*}
\Xi(2)\sum_{n=2}^{\infty} |a_{n}| \leq \sum_{n=2}^{\infty}\Xi(n) \,|a_{n}| \leq (A-B)(1-\gamma),
\end{align*}
that is
\begin{align*}
\sum_{n=2}^{\infty} |a_{n}| \leq \frac{(A-B)(1-\gamma)}{\Xi(2)}.
\end{align*}
Therefore, we have
\begin{align*}
|F_\mu(z)|\leq |z|+ |z|^{2 \mu} \sum_{n=2}^{\infty} |a_{n}|,
\end{align*}
\begin{align*}
|F_\mu(z)|\leq |z|+ \frac{(A-B)(1-\gamma)}{\big[(1-B) \left(( 2\mu
)^{k} \vartheta_{2} - (2\mu )^{m} \lambda_{2}\right)
 + (A-B)(1-\gamma)(2\mu)^{m} \vartheta_{2}\big]} |z|^{ 2 \mu }.
\end{align*}
Similarly, we have
\begin{align*}
|F_\mu(z)| &\geq |z|- \sum_{n=2}^{\infty} |a_{n}| |z|^{\mu n} \geq
|z| - |z|^{2\mu}\sum_{n=2}^{\infty} |a_{n}| \nonumber
\\
&\geq |z|- \frac{(A-B)(1-\gamma)}{\big[(1-B) \left(( 2\mu )^{k} \vartheta_{2} - (2\mu )^{m} \lambda_{2}\right)
 + (A-B)(1-\gamma)(2\mu)^{m} \vartheta_{2}\big]} |z|^{ 2 \mu }.
\end{align*}
The result is sharp for the function
\begin{align}\label{F*}
F_\mu(z)= |z|- \frac{(A-B)(1-\gamma)}{\big[(1-B) \left(( 2\mu )^{k}
\vartheta_{2} - (2\mu )^{m} \lambda_{2}\right)
 + (A-B)(1-\gamma)(2\mu)^{m} \vartheta_{2}\big]} |z|^{ 2 \mu }.
\end{align}
This completes the proof of Theorem  \ref{Th:7}.
 \end{proof}
}\end{theorem}

\begin{theorem}\label{The8}\rm{
Let $F_\mu(z)$ defined by \eqref{-} be in the class $\widetilde{
\mathcal{E}}_{k,m}(\Phi_\mu,\Psi_\mu,A,B,\mu,\gamma)$, then
\begin{align}
|F^{\prime}_\mu(z)| &\geq 1- \frac{2 \mu (A-B)(1-\gamma)}{\big[(1-B)
\left(( 2\mu )^{k} \vartheta_{2} - (2\mu )^{m} \lambda_{2}\right)
 + (A-B)(1-\gamma)(2\mu)^{m} \vartheta_{2}\big]} |z|^{ 2 \mu-1 } \\
&\leq  1 +  \frac{2 \mu (A-B)(1-\gamma)}{\big[(1-B) \left(( 2\mu )^{k} \vartheta_{2} - (2\mu )^{m} \lambda_{2}\right)
 + (A-B)(1-\gamma)(2\mu)^{m} \vartheta_{2}\big]} |z|^{ 2 \mu-1 }
\end{align}
This result is sharp.
\begin{proof}
Similarly $ \Xi(n)/n $ is an increasing function of $ n (n\geq 1)$,
\begin{align}
 \frac{\Xi(2)}{2\mu }\sum_{n=2}^{\infty}  (\mu n) \, |a_{n}| & \leq
  \sum_{n=2}^{\infty} \frac{\Xi(n)}{(\mu n)} (\mu n)\, |a_{n}|
   = \sum_{n=2}^{\infty} \Xi(n) |a_{n}| \leq (1-\gamma)(A-B),
\end{align}
that is
\begin{align}
\sum_{n=2}^{\infty} (\mu n) |a_{n}| \leq \frac{2 \mu ( 1-\gamma)(A-B)}{\Xi(2)}.
\end{align}
Then, we get
\begin{align}
|F^{\prime}_\mu(z)| & \leq 1+ |z|^{2\mu-1} \sum_{n=2}^{\infty} (\mu n) |a_{n}| \nonumber \\
& \leq
 1+ \frac{2 \mu\, (A-B)(1-\gamma)} {\big[(1-B) \left(( 2\mu )^{k} \vartheta_{2} - (2\mu )^{m} \lambda_{2}\right)
 + (A-B)(1-\gamma)(2\mu)^{m} \vartheta_{2}\big]} |z|^{ 2 \mu-1 }.
\end{align}
and similarly
\begin{align}
|F^{\prime}_\mu(z)| & \geq 1 - |z|^{2\mu-1} \sum_{n=2}^{\infty} (\mu n) |a_{n}| \nonumber \\
& \geq
 1- \frac{2 \mu\, (A-B)(1-\gamma)}{\big[(1-B) \left(( 2\mu )^{k} \vartheta_{2} - (2\mu )^{m} \lambda_{2}\right)
 + (A-B)(1-\gamma)(2\mu)^{m} \vartheta_{2}\big]} |z|^{ 2 \mu-1 }.
\end{align}
It is clear that the assertions of Theorem \ref{The8} are sharp for
the function $F_\mu(z)$ given by \eqref{F*}. This complete the proof
of Theorem \ref{The8}.
\end{proof}
}\end{theorem}


\subsection{The radii subclasses of class $ \widetilde{ \mathcal{E}}_{k,m}(\Phi_\mu,\Psi_\mu,A,B,\mu,\gamma) $ }
In this section radii of bounded turning, convexity and starlikeness
for functions $ F_\mu(z) \in \widetilde{ \mathcal{E}}_{k,m}
(\Phi_\mu,\Psi_\mu,A,B,\mu,\gamma) $ are studied.

\begin{theorem}\rm{
Let $ F_\mu(z)$ given by \eqref{-} be in the class  $\widetilde{
\mathcal{E}}_{k,m}(\Phi_\mu,\Psi_\mu,A,B,\mu,\gamma) $, then

\begin{itemize}
\item[i-] $ F_\mu(z)$ is starlike of order $\psi \,( 0 \leq \psi < 1)$ in $ |z| < r_{1}$, where
\begin{flalign}\label{r1}
r_{1}= \inf_{n\geq 2} \lbrace & \frac{\big[(1-B) \left((\mu n)^{k} \vartheta_{n} - (\mu n)^{m} \lambda_{n}\right)
 + (A-B)(1-\gamma)(\mu n)^{m} \vartheta_{n}\big]}  {(A-B)(1-\gamma)} \nonumber
 \\
 &  \qquad \qquad \qquad \qquad \qquad \qquad  \times \left( \frac{1-\psi}{\mu n - \psi}\right)\rbrace^{1/(\mu n -1)}.
 \end{flalign}
 \item[ii-] $F_\mu(z)$ is convex of order $ \psi( 0 \leq \psi < 1)$ in $|z| < r_{2}$, where
 \begin{align}\label{r2}
r_{2}= \inf_{n\geq 2} \lbrace&  \frac{\big[(1-B) \left((\mu n)^{k} \vartheta_{n} - (\mu n)^{m} \lambda_{n}\right)
 + (A-B)(1-\gamma)(\mu n)^{m} \vartheta_{n}\big]}  {(A-B)(1-\gamma)}\nonumber
 \\
 & \qquad \qquad \qquad \qquad \qquad \qquad  \times \left( \frac{1-\psi}{\mu n(\mu n - \psi)}\right) \rbrace^{1/(\mu n -1)},
 \end{align}
 \item[iii-]$F_\mu(z)$ is close to convex of order $ \psi( 0 \leq \psi < 1)$ in $|z| < r_{2}$, where
 \begin{align}\label{r3}
r_{3}= \inf_{n\geq 2} \lbrace&  \frac{\big[(1-B) \left((\mu n)^{k} \vartheta_{n} - (\mu n)^{m} \lambda_{n}\right)
 + (A-B)(1-\gamma)(\mu n)^{m} \vartheta_{n}\big]}  {(A-B)(1-\gamma)}\nonumber
 \\
 & \qquad \qquad \qquad \qquad \qquad \qquad  \times \left( \frac{1-\psi}{\mu n }\right) \rbrace^{1/(\mu n -1)},
 \end{align}
 \end{itemize}
 Each the results is sharp for the function $ F_\mu(z)$ given by \eqref{fun}

\begin{proof}
It is sufficient show that
\begin{align} \label{Eq5}
\left | \frac{z F_\mu^{\prime}(z)}{F_\mu(z)}\right| \leq 1-\psi,
\quad \text{for} \quad |z|< r_{1},
\end{align}
where $ r_{1}$ is defined by \eqref{r1}. Further, we find from \eqref{-} that
\begin{align}
\left| \frac{zF_\mu^{\prime}(z)}{F_\mu(z)} -1 \right|\leq
\frac{\sum_{n=2}^{\infty} ( \mu n - 1) a_{n} |z|^{\mu n -1} }{1-
\sum_{n=2}^{\infty} a_{n} |z|^{\mu n-1}}.
\end{align}
Therefore, we satisfy \eqref{Eq5} if and only if
\begin{align}\label{Iq1}
\sum_{n=2}^{\infty} \frac{(\mu - \psi) a_{n} |z|^{\mu n-1}}{(1-\psi)} \leq 1.
\end{align}
Nevertheless, from Theorem \ref{Th1}, inequality
\eqref{Iq1} it will be true that if
\begin{align}
\frac{(\mu - \psi)  |z|^{\mu n-1}}{(1-\psi)} \leq \frac{\big[(1-B) \left((\mu n)^{k} \vartheta_{n} - (\mu n)^{m} \lambda_{n}\right)
 + (A-B)(1-\gamma)(\mu n)^{m} \vartheta_{n}\big]}  {(A-B)(1-\gamma)}\nonumber
\end{align}
this is, if
\begin{align}
|z| \leq \lbrace \frac{\big[(1-B) \left((\mu n)^{k} \vartheta_{n} - (\mu n)^{m} \lambda_{n}\right)
 + (A-B)(1-\gamma)(\mu n)^{m} \vartheta_{n}\big]}  {(A-B)(1-\gamma)}\nonumber  \nonumber  \\
  \times \left( \frac{1-\psi}{\mu n - \psi}\right)\rbrace^{1/(\mu n -1)}
\end{align}
Or equivalent to
\begin{align}
r_{1}= \inf_{n\geq 2}  \lbrace \frac{\big[(1-B) \left((\mu n)^{k} \vartheta_{n} - (\mu n)^{m} \lambda_{n}\right)
 + (A-B)(1-\gamma)(\mu n)^{m} \vartheta_{n}\big]}  {(A-B)(1-\gamma)}\nonumber \nonumber  \\
  \times \left( \frac{1-\psi}{\mu n - \psi}\right)\rbrace^{1/(\mu n -1)}.
\end{align}
This completes the  proof of \eqref{r1}.
To prove \eqref{r2} and \eqref{r3}; respectively  it is sufficient to show that
\begin{align}
\left| 1 + \frac{z F_\mu^{\prime \prime}(z)}{F_\mu^{\prime}(z)} -1
\right| \leq 1- \psi \quad (|z| \leq r_{2}; \quad 0 \leq \psi < 1)
\end{align}
and
\begin{align}
\left|   F_\mu^{\prime}(z)-1 \right| \leq 1- \psi \quad (|z| \leq
r_{3}; \quad 0 \leq \psi < 1).
\end{align}
\end{proof}
}\end{theorem}


\subsection{Extreme point}
\begin{theorem}\label{Th10}\rm{
Let $ F_{1}(z) =z,$ and
\begin{align}
F_{\mu n}(z)= z- \frac{(A-B)(1-\gamma)}{\big[(1-B) \left(( \mu n )^{k} \vartheta_{n} - (\mu n )^{m} \lambda_{n}\right)
 + (A-B)(1-\gamma)(\mu n)^{m} \vartheta_{n}\big]}   z^{ \mu n }
\end{align}

Then $ F_\mu(z) \in \widetilde{
\mathcal{E}}_{k,m}(\Phi,\Psi,A,B,\mu,\gamma)$ if and only if it can
be expressed in the following form:
$$ F(z)= \sum_{n=1}^{\infty} \eta_{n} F_{\mu n}(z),$$
where
$$ \eta_{n} \geq 0, \quad \text{and} \quad \sum_{n=1}^{\infty}\eta_{n}= 1.$$

\begin{proof}
Assuming that,
\begin{align}
F_\mu(z)  & = \sum_{n=1}^{\infty} \eta_{n} F_{\mu n}(z)\\
&= z -\sum_{n=1}^{\infty} \eta_{n} \frac{(A-B)(1-\gamma)}{\big[(1-B)
\left(( \mu n )^{k} \vartheta_{n} - (\mu n )^{m} \lambda_{n}\right)
 + (A-B)(1-\gamma)(\mu n)^{m} \vartheta_{n}\big]}   z^{ \mu n }
\end{align}
Then, from Theorem \ref{Th1}, we obtain
\begin{align}
  \sum_{n=2}^{\infty}& \bigg({\big[(1-B)\left((\mu n)^{k} \vartheta_{n}-(\mu n)^{m} \lambda_{n}\right)+(A-B)(1-\gamma)(\mu n)^{m} \vartheta_{n}\big]}
  \nonumber \\
  & \qquad \times \frac{(A-B)(1-\gamma)}{\big[(1-B) \left(( \mu n )^{k} \vartheta_{n} - (\mu n )^{m} \lambda_{n}\right)
 + (A-B)(1-\gamma)(\mu n)^{m} \vartheta_{n}\big]}\bigg)\nonumber \\
 & = (A-B)(1-\gamma)\sum_{n=2}^{\infty} \eta_{n}= (A-B)(1-\gamma)(1-\eta_{1}) \leq (A-B)(1-\gamma).\nonumber
 \end{align}
 Therefore, in view of Theorem\ref{Th1}, we find that $F(z) \in \widetilde{ \mathcal{E}}_{k,m}(\Phi,\Psi,A,B,\mu,\gamma)$.
 Conversely, let us suppose that $ F(z) \in  \widetilde{ \mathcal{E}}_{k,m}(\Phi,\Psi,A,B,\mu,\gamma)$, then, since
\begin{align}
a_{n}\leq \frac{(A-B)(1-\gamma)}{\big[(1-B) \left(( \mu n )^{k} \vartheta_{n} - (\mu n )^{m} \lambda_{n}\right)
 + (A-B)(1-\gamma)(\mu n)^{m} \vartheta_{n}\big]}
\end{align}
by setting
\begin{align*}
 \eta_{n}= \frac {\big[(1-B)\left((\mu n)^{k} \vartheta_{n}-(\mu n)^{m}
 \lambda_{n}\right)+(A-B)(1-\gamma)(\mu n)^{m} \vartheta_{n}\big]} {(A-B)(1-\gamma)}a_{n},\, n =\lbrace2,3,..\rbrace
\end{align*}
 and
 $$ \eta_{1}= 1 - \sum_{n=2}^{\infty} \eta_{n}.$$
Therefore, we have
\begin{align*}
F(z)= \sum_{n=1}^{\infty} \eta_{n} F_{\mu n}(z).
\end{align*}
By this, we complete the proof of Theorem \ref{Th10}.
\end{proof}
}\end{theorem}
\begin{corollary}\rm{
The extreme point of the class $\widetilde{ \mathcal{E}}_{k,m}(\Phi,\Psi,A,B,\mu,\gamma) $ are given by
\begin{align*}
F_{\mu n}(z)=z-\frac{(A-B)(1-\gamma)}{\big[(1-B) \left(( \mu n )^{k} \vartheta_{n} - (\mu n )^{m} \lambda_{n}\right)
 + (A-B)(1-\gamma)(\mu n)^{m} \vartheta_{n}\big]} z^{\mu n}.
\end{align*}
}\end{corollary}

\subsection{Integral Means Inequality}
In this section, we consider some result due to Littlewood subordination (see\cite{z2}).

\begin{lemma}\label{Lemma 3}\rm{
If the functions $ f$ and $g$ are analytic in open unit disk $ \mathbb{U}$ with
\begin{align}
f(z) \prec g(z) \quad ( z \in \mathbb{U}),
\end{align}
then, for $ q>0$ and $ z= e^{i\theta}\quad ( 0 < r< 1),$
\begin{align}\label{55}
\int_{0}^{2\pi}|f(z)|^{q} d\theta \leq \int_{0}^{2\pi}|g(z)|^{q} d\theta.
\end{align}
}\end{lemma} \noindent
Now, let use Lemma\ref{Lemma 3} to prove the following Theorem.

\begin{theorem} \rm{ Assume that $ F_\mu(z) \in \widetilde{ \mathcal{E}}_{k,m}(\Phi_\mu,\Psi_\mu,A,B,\mu,\gamma)
$, $  q> 0, -1 \leq B< A \leq 1$, $ k; m \in \mathbb{N}_{0}$, and  $
F_{2\mu}(z)$ is defined by
 \begin{align}
F_{2\mu}(z)= z- \frac{(A-B)(1-\gamma)}{\big[(1-B) \left(( 2\mu  )^{k} \vartheta_{2} - (2\mu  )^{m} \lambda_{2}\right)
 + (A-B)(1-\gamma)(2\mu)^{m} \vartheta_{2}\big]} z^{2\mu}.
\end{align}
Then $ z= r e^{i \theta} ( 0 < r< 1)$, we obtain
$$ \int_{0}^{2 \pi} |F_\mu(z)|^{q}d\theta \leq \int_{0}^{2\pi} |F_{2\mu}(z)|^{q} d\theta.$$
\begin{proof}
For $F_\mu(z)=z- \sum_{n=2}^{\infty}a_{n} z^{\mu n} \quad(a_{n} \geq
0)$ and by \eqref{55} is equivalent to proof,
\begin{align}
&\int_{0}^{2\pi}\left| 1- \sum_{n=2}^{\infty} a_{n}z^{\mu n-1} \right|^{q} d \theta  \nonumber
 \\
&
\leq \int_{0}^{2\pi}\left| 1- \frac{(A-B)(1-\gamma)}{\big[(1-B) \left(( 2\mu  )^{k} \vartheta_{2} - (2\mu  )^{m} \lambda_{2}\right)
 + (A-B)(1-\gamma)(2\mu)^{m} \vartheta_{2}\big]} z^{2\mu-1}\right|^{q} d \theta
\end{align}
By applying Lemma \ref{Lemma 3}, it would suffice to show that
\begin{align*}
1-&\sum_{n=2}^{\infty}a_{n}z^{\mu n-1}\\
&
\prec
1-\frac{(A-B)(1-\gamma)}{\big[(1-B) \left(( 2\mu  )^{k} \vartheta_{2} - (2\mu  )^{m} \lambda_{2}\right)
 + (A-B)(1-\gamma)(2\mu)^{m} \vartheta_{2}\big]} z^{2\mu-1}.
\end{align*}
By putting
\begin{align*}
1- & \sum_{n=2}^{\infty} a_{n} z^{\mu n-1} \\
&
=
1-\frac{(A-B)(1-\gamma)}{\big[(1-B) \left(( 2\mu  )^{k} \vartheta_{2} - (2\mu  )^{m} \lambda_{2}\right)
 + (A-B)(1-\gamma)(2\mu)^{m} \vartheta_{2}\big]} \Theta(z).
\end{align*}
and by using Theorem\ref{Th1}, we have
\begin{align*}
 & \left| \sum_{n=2}^{\infty} \frac{\big[(1-B) \left(( 2\mu  )^{k} \vartheta_{2} - (2\mu  )^{m} \lambda_{2}\right)
 + (A-B)(1-\gamma)(2\mu)^{m} \vartheta_{2}\big]}{(A-B)(1-\gamma)}  a_{n} z^{2\mu-1}  \right|
\nonumber \\
 & \leq |z|^{2\mu-1} \sum_{n=2}^{\infty} \frac{\big[(1-B) \left(( 2\mu  )^{k} \vartheta_{2} - (2\mu  )^{m} \lambda_{2}\right)
 + (A-B)(1-\gamma)(2\mu)^{m} \vartheta_{2}\big]}{(A-B)(1-\gamma)}  a_{n}
 \nonumber\\
 & \leq |z|^{2\mu-1} \sum_{n=2}^{\infty} \frac{\big[(1-B) \left(( 2\mu  )^{k} \vartheta_{2} - (2\mu  )^{m} \lambda_{2}\right)
 + (A-B)(1-\gamma)(2\mu)^{m} \vartheta_{2}\big]}{(A-B)(1-\gamma)}  a_{n}  \nonumber\\
 & \leq |z| \leq  1.
\end{align*}
This complete the proof of this Theorem.
\end{proof}
}\end{theorem}


\section{Distortion Applications }
In this section, we prove some distortion theorems involving certain fractional calculus (Srivastava and Owa)
operators for function $ F(z)$ belonging to class $\widetilde{ \mathcal{E}}_{k,m}(\Phi,\Psi,A,B,\mu,\gamma) $.
 Now let recall the following definitions:

\begin{definition}\label{d1} \rm{ The fractional derivative of order $ \delta$ is defined, for a function $ f(z)$ by
\begin{align}
\mathcal{D}_{z}^{\delta} f(z) = \frac{1}{\Gamma(1-\delta)}\frac{d}{dz}\int_{0}^{z} (z-\zeta)^{-\delta} f(\zeta) d\zeta
\quad ( 0 \leq \delta < 1),
\end{align}
where the function $ f(z)$ is constrained and the multiplicity of $ (z-\zeta)^{-\delta}$ is removed  by requiring
$ \log(z-\zeta)$ to be real when $ z-\zeta>0$(see\cite{z7},\cite{z10}).
}\end{definition}
   \begin{definition}\label{d2}\rm{ The fractional integral of order $ \delta$ is defined, for a function $ f(z)$ by
  \begin{align}
\mathcal{I}_{z}^{\delta} f(z) = \frac{1}{\Gamma(\delta)}\frac{d}{dz}\int_{0}^{z} (z-\zeta)^{\delta-1} f(\zeta) d\zeta \quad ( \delta >0),
\end{align}
where the function $ f(z)$ is  analytic in a simply connected domain of the complex $ z$-plane containing the origin,
 and the multiplicity of  $ (z-\zeta)^{-\delta}$ is removed by suggesting $ \log(z-\zeta)$ to be real
  when $z-\zeta>0$ (see\cite{z7},\cite{z10}).
}\end{definition}
\begin{definition}\label{d3}\rm{(see\cite{z7})
Under the hypotheses of Definition \ref{d1}, the fractional derivative of order $ \delta$ is defined for a function $ f(z)$, by
\begin{align}
\mathcal{D}_{z}^{\upsilon+\delta} = \frac{d^{\upsilon}}{dz^{\upsilon}}\lbrace\mathcal{D}_{z}^{\delta}\rbrace \quad
( 0 \leq \delta < 1; \upsilon\in \lbrace0,1,2,..\rbrace).
\end{align}
}\end{definition}
Using Definitions \ref{d1}, \ref{d2} and \ref{d3}, we have
\begin{lemma}\rm{
For the function $z^{\mu n}$, $(\mu \geq 1; n \in \mathbb{N}_{0})$, $z \in \mathbb{U}$, we have

\begin{align}
\mathcal{D}_{z}^{\delta}\lbrace z^{\mu n}\rbrace = \frac{\Gamma(\mu n + 1)}{\Gamma(\mu n + 1 - \delta)}z^{\mu n - \delta}\quad  (0 \leq \delta < 1)
\end{align}
and
\begin{align}
\mathcal{I}_{z}^{\delta}\lbrace z^{\mu n}\rbrace = \frac{\Gamma(\mu n + 1)}
{\Gamma(\mu n + 1 + \delta)}z^{\mu n+ \delta}\quad  (0 \leq \delta < 1).
\end{align}
}\end{lemma}
\begin{theorem}\label{33}\rm{
Let the function $ F_\mu(z)$ given by \eqref{-} be in the class
$\widetilde{ \mathcal{E}}_{k,m}(\Phi_\mu,\Psi_\mu,A,B,\mu,\gamma)$.
Then
\begin{align*}
&|\mathcal{I}_{z}^{\delta}F_\mu(z)| \geq \frac{|z|^{\mu+\delta}}{\Gamma(2+\delta)}\\
& \times \left(1-
\frac{(1)_{2\mu-1}(A-B)(1-\gamma)}{(2-\delta)_{2\mu-1}\lbrace
(1-B)[(2\mu )^{k-1} \vartheta_{2} - (2\mu )^{m-1} \lambda_{2}]  +
(A-B)(1-\gamma)(2\mu )^{k-1} \vartheta_{2}\rbrace}|z|^{\mu} \right)
\end{align*}
and
\begin{align*}
&\mathcal{I}_{z}^{\delta}F_\mu(z)\leq  \frac{|z|^{\mu+\delta}}{\Gamma(2+\delta)}\\
& \times \left(1+
\frac{(1)_{2\mu-1}(A-B)(1-\gamma)}{(2-\delta)_{2\mu-1}\lbrace
 (1-B)[(2\mu )^{k-1} \vartheta_{2} - (2\mu )^{m-1} \lambda_{2}]  + (A-B)(1-\gamma)(2\mu )^{k-1} \vartheta_{2}\rbrace}|z|^{\mu} \right).
\end{align*}
The results are sharp.
\begin{proof}
Let
\begin{align}
\mathcal{F}(z) & = \Gamma(2+\delta) z^{-\delta} \mathcal{I}_{z}^{\delta} F_\mu(z) \nonumber \\
&= z- \sum_{n=2}^{\infty} \frac{\Gamma(\mu n+1) \Gamma(2-\delta)}{\Gamma(\mu n+1 +\delta)} a_{n} z^{\mu n},\\
&=z- \sum_{n=2}^{\infty} \Upsilon(n) a_{n} z^{\mu n}
\end{align}
where
\begin{align}
\Upsilon(n)= \frac{\Gamma(\mu n+1) \Gamma(2+\delta)}{\Gamma(\mu n+1 +\delta)}\quad ( n =2,3,...).
\end{align}
It is clear that $ \Upsilon(n)$ is a decreasing function of $ n$, we can write as
\begin{align}\label{85}
0 < \Upsilon(n) \leq \Upsilon(2) = \frac{2\mu (1)_{2\mu-1}}{(2+\delta)_{2\mu-1}}.
\end{align}
On the other hand, in view of Theorem \ref{Th2}, we obtain
\begin{align*}
& \big[(1-B)\left((2\mu )^{k} \vartheta_{2} - (2\mu )^{m} \lambda_{2}\right)  + (A-B)(1-\gamma)(2\mu )^{k} \vartheta_{2}\big] \sum_{n=2}^{\infty}
 a_{n} \\
 & \leq
 \sum_{n=2}^{\infty} \big[(1-B)\left((\mu n )^{k} \vartheta_{n} - (\mu n )^{m} \lambda_{n}\right)  + (A-B)(1-\gamma)(\mu n )^{m} \vartheta_{2}\big]
 a_{n} \\
 &
 \leq (A-B)(1-\gamma).
\end{align*}
Then
\begin{align}\label{86}
\sum_{n=2}^{\infty} a_{n} \leq \frac{(A-B)(1-\gamma)}{\big[(1-B)\left((2\mu )^{k} \vartheta_{2} - (2\mu )^{m} \lambda_{2}\right)  + (A-B)(1-\gamma)(2\mu )^{m} \vartheta_{2}\big]}
\end{align}
Thus, by using \eqref{85} and \eqref{86}, we see that
\begin{align}
|\mathcal{F}(z)| & \geq |z|- \Upsilon(2) |z|^{2\mu} \sum_{n=2}^{\infty} a_{n}
\\
& \geq |z|- \frac{(1)_{2\mu-1}(A-B)(1-\gamma)}{(2-\delta)_{2\mu-1}\big[(1-B)\left((2\mu )^{k-1} \vartheta_{2} - (2\mu )^{m-1} \lambda_{2}\right)  + (A-B)(1-\gamma)(2\mu )^{m-1} \vartheta_{2}\big]}|z|^{2\mu}
\end{align}
Similarly,
\begin{align}
|\mathcal{F}(z)| & \leq |z|+ \Upsilon(2) |z|^{2\mu} \sum_{n=2}^{\infty} a_{n}
\\
& \leq |z| + \frac{(1)_{2\mu-1}(A-B)(1-\gamma)}{(2-\delta)_{2\mu-1}\big[(1-B)\left((2\mu )^{k-1} \vartheta_{2} -
(2\mu )^{m-1} \lambda_{2}\right)  + (A-B)(1-\gamma)(2\mu )^{m-1} \vartheta_{2}\big]}|z|^{2\mu}
\end{align}
From, this prove Theorem \ref{33}.
The equalities are attained for the function $ F(z)$ given by
\begin{align*}
&|\mathcal{I}_{z}^{\delta}F_\mu(z)|= \frac{z^{\mu+\delta}}{\Gamma(2+\delta)}\\
& \times \left(1-
\frac{(1)_{2\mu-1}(A-B)(1-\gamma)}{(2-\delta)_{2\mu-1}\big[(1-B)\left((2\mu
)^{k-1} \vartheta_{2} - (2\mu )^{m-1} \lambda_{2}\right)  +
(A-B)(1-\gamma)(2\mu )^{m-1} \vartheta_{2}\big]}|z|^{2\mu}\right)
\end{align*}
Then the result are sharp and the proof of Theorem \ref{33} is completed.
\end{proof}
}\end{theorem}
\begin{corollary}\rm{
Under the hypothesis of Theorem \ref{33}, $
\mathcal{I}_{z}^{\delta}F_\mu(z)$ is included in a disk with its
center at the origin and radius $ \mathcal{R}_{1}$ given by
\begin{align}
&\mathcal{R}_{1}= \frac{1}{\Gamma(2+\delta)}\\
& \times\left(1-   \frac{(1)_{2\mu-1}(A-B)(1-\gamma)}{(2-\delta)_{2\mu-1}\big[(1-B)\left((2\mu )^{k-1} \vartheta_{2} -
 (2\mu )^{m-1} \lambda_{2}\right)  + (A-B)(1-\gamma)(2\mu )^{m-1} \vartheta_{2}\big]}\right)
\end{align}
}\end{corollary}
\begin{theorem}\label{44}\rm{
Let the function $ F_\mu(z)$ given by \eqref{-} be in the class $
\mathcal{E}_{k,m}(\Phi_\mu,\Psi_\mu,A,B,\gamma,\mu) $, then
\begin{align*}
|\mathcal{D}_{z}^{\delta}&F_\mu(z)|\geq \frac{|z|^{\mu-\delta}}{\Gamma(2-\delta)} \\
&\times \left(1-
\frac{(1)_{2\mu-1}(A-B)(1-\gamma)}{(2-\delta)_{2\mu-1}\big[(1-B)\left((2\mu
)^{k-1} \vartheta_{2} - (2\mu )^{m-1} \lambda_{2}\right)  +
(A-B)(1-\gamma)(2\mu )^{m-1} \vartheta_{2}\big]}|z|^{2\mu}\right)
\end{align*}
and
\begin{align*}
|\mathcal{D}_{z}^{\delta}&F_\mu(z)| \leq  \frac{|z|^{\mu-\delta}}{\Gamma(2-\delta)} \\
&\times \left(1-
\frac{(1)_{2\mu-1}(A-B)(1-\gamma)}{(2-\delta)_{2\mu-1}\big[(1-B)\left((2\mu
)^{k-1} \vartheta_{2} - (2\mu )^{m-1} \lambda_{2}\right)  +
(A-B)(1-\gamma)(2\mu )^{m-1} \vartheta_{2}\big]}|z|^{2\mu}\right).
\end{align*}
Each f these results is sharp.
\begin{proof}
Let
\begin{align*}
\mathcal{G}(z) & = \Gamma(2-\delta) z^{\delta} \mathcal{I}_{z}^{\delta} F_\mu(z) \nonumber \\
&= z- \sum_{n=2}^{\infty} \frac{\Gamma(\mu n) \Gamma(2-\delta)}{\Gamma(\mu n+1 +\delta)} (\mu n) a_{n} z^{\mu n},\\
&=z- \sum_{n=2}^{\infty} \Omega(n)(\mu n) a_{n} z^{\mu n}
\end{align*}
where
\begin{align}
\Omega(n)= \frac{\Gamma(\mu n) \Gamma(2-\delta)}{\Gamma(\mu n+1 -\delta)}\quad ( n =2,3,...).
\end{align}
It is clear that $ \Omega(n)$ is a decreasing function of $ n$, we can write as
\begin{align}\label{85}
0 < \Omega(n) \leq \Omega(2) = \frac{ (1)_{2\mu-1}}{(2-\delta)_{2\mu-1}}.
\end{align}
On the other hand, in view of Theorem \ref{Th2}, we obtain
\begin{align*}
\big[&(1-B) \left((2\mu )^{k} \vartheta_{2} - (2\mu )^{m}
\lambda_{2}\right)  + (A-B)(1-\gamma)(2\mu )^{m} \vartheta_{2}\big] \sum_{n=2}^{\infty} (\mu n) a_{n} \\
 & \leq
 \sum_{n=2}^{\infty} \big[(1-B)\left((\mu n)^{k} \vartheta_{n} - (\mu n)^{m} \lambda_{n}\right)
 + (A-B)(1-\gamma)(\mu n)^{m}  \vartheta_{n}\big] a_{n} \\
 &
 \leq (A-B)(1-\gamma).
\end{align*}
Then
\begin{align}\label{86}
\sum_{n-2}^{\infty}(\mu n)a_{n} \leq   \frac{(A-B)(1-\gamma)}{\big[(1-B)\left((2\mu )
^{k-1} \vartheta_{2} - (2\mu )^{m-1} \lambda_{2}\right)  + (A-B)(1-\gamma)(2\mu )^{m-1} \vartheta_{2}\big]}
\end{align}
Thus, by sing \eqref{85} and \eqref{86}, we see that
\begin{align}
|\mathcal{G}(z)| & \geq |z|- \Omega(2) |z|^{2\mu} \sum_{n=2}^{\infty} (\mu n)a_{n}
\\
& \geq |z|- \frac{(1)_{2\mu-1}(A-B)(1-\gamma)}{(2-\delta)_{2\mu-1}\big[(1-B)
\left((2\mu )^{k-1} \vartheta_{2} - (2\mu )^{m-1} \lambda_{2}\right)
+ (A-B)(1-\gamma)(2\mu )^{m-1} \vartheta_{2}\big]}|z|^{2\mu}
\end{align}
Similarly,
\begin{align}
|\mathcal{G}(z)| & \leq |z|+ \Omega(2) |z|^{2\mu} \sum_{n=2}^{\infty}(\mu n) a_{n}
\\
& \leq |z| + \frac{(1)_{2\mu-1}(A-B)(1-\gamma)}{(2-\delta)_{2\mu-1}\big[(1-B)\left((2\mu )^{k-1} \vartheta_{2} - (2\mu )^{m-1}
 \lambda_{2}\right)  + (A-B)(1-\gamma)(2\mu )^{m-1} \vartheta_{2}\big]}|z|^{2\mu}
\end{align}
From, this prove Theorem \ref{44}. The equalities are attained for
the function $ F_\mu(z)$ given by
\begin{align*}
|\mathcal{D}_{z}^{\delta}&F_\mu(z)|= \frac{z^{2\mu-\delta}}{\Gamma(2-\delta)} \\
&\times \left(1-
\frac{(1)_{2\mu-1}(A-B)(1-\gamma)}{(2-\delta)_{2\mu-1}\big[(1-B)
\left((2\mu )^{k-1} \vartheta_{2} - (2\mu )^{m-1} \lambda_{2}\right)
+ (A-B)(1-\gamma)(2\mu )^{m-1} \vartheta_{2}\big]}|z|^{2\mu} \right)
\end{align*}
Then the result are sharp and the proof of Theorem \ref{33} is completed.
\end{proof}
}\end{theorem}
\begin{corollary}\rm{
Under the hypothesis of Theorem \ref{44},
 $ \mathcal{D}_{z}^{\delta}F_\mu(z)$ is included in a disk with its center at the origin and radius $ \mathcal{R}_{2}$ given by
\begin{align}
\mathcal{R}_{1}&= \frac{1}{\Gamma(2-\delta)}\\
& \times \left(1-   \frac{(1)_{2\mu-1}(A-B)(1-\gamma)}
{(2-\delta)_{2\mu-1}\big[(1-B)\left((2\mu )^{k-1} \vartheta_{2} - (2\mu )^{m-1}
\lambda_{2}\right)  + (A-B)(1-\gamma)(2\mu )^{m-1} \vartheta_{2}\big]}\right)
\end{align}
}\end{corollary}

  \section{Conclusion} In unit disk, we derived a new class of fractional power $ \mathcal{A}_{\mu}$
  and consider this class to define a generalized of many
   differential operator, also we employed this operator to define
   a new subclasses in open unit disk. Further, we studied some
    Characteristic properties and applications including certain
     fractional calculus operators.

  \end{document}